\documentclass[a4paper]{article}
\usepackage[english]{babel}
\usepackage[utf8]{inputenc}
\usepackage{amsmath,amsthm}
\usepackage{amsfonts,amssymb}

\usepackage{tikz,tikz-cd} 

\usepackage{xifthen} 

\newcommand{\toposdefaut}{0}
\newcommand{\topos}[1][\toposdefaut]{ 
\ifthenelse{\equal{#1}{0}}{ \mathcal{T} }
{
\ifthenelse{\equal{#1}{1}}{ \mathcal{E} }{ #1 }
}
}

\newcommand{\sh}{\textsf{Sh}}

\newcommand{\spec}{\text{Spec }}

\newcommand{\Ecal}{\mathcal{E}} 
 
\newcommand{\Tcal}{\mathcal{T}}

\newcommand{\Ocal}{\mathcal{O}}

\newcommand{\Scal}{\mathcal{S}} 
 
\newcommand{\Fcal}{\mathcal{F}}

\newcommand{\Lcal}{\mathcal{L}}

\newcommand{\Ccal}{\mathcal{C}}

\newcommand{\Dgo}{\mathfrak{D}}

\usepackage{titlesec}
\titleformat{\subsubsection}[runin]{\normalfont}{\thesubsubsection}{0pt}{}[.]

\renewcommand{\thesubsubsection}{\arabic{section}.\arabic{subsubsection}}
\csname @addtoreset\endcsname{subsubsection}{section}

\newcommand{\block}[1]
{

\par \subsubsection{} #1

\bigskip}

\newcommand{\blockn}[1]{\par #1 \bigskip}

\newcommand{\Th}[1]
	{
	\bigskip	
	\textbf{Theorem : }{\itshape #1}
		
	\bigskip
	}

\newcommand{\Prop}[1]
	{

	\bigskip
	
	\textbf{Proposition : }{\itshape #1}
		
	\bigskip
	
	}

\newcommand{\Cor}[1]
	{

	\bigskip
	
	\textbf{Corollary : }{\itshape #1}	
		
	\bigskip

	}

\newcommand{\Lem}[1]
	{

	\bigskip
	
	\textbf{Lemma : }{\itshape #1}
		
	\bigskip
	
	}

\newcommand{\Def}[1]
	{
	
	\bigskip
	
	\textbf{Definition : }{\itshape #1}
	
	\bigskip
	
	}

\newcommand{\Dem}[1]{
	
	\smallskip
	
	\textbf{Proof : } \par
	 {#1} $\square$
	 
	 \bigskip
}

\begin{document}

\pagestyle{plain}
\title{On toposes generated by cardinal finite objects}
\author{Simon Henry}
\date{}

\maketitle

\begin{abstract}
We give a characterizations of toposes which admit a generating set of objects which are internally cardinal finite (i.e. Kuratowski finite and decidable) in terms of ``topological" conditions. The central result is that, constructively, a hyperconnected separated locally decidable topos admit a generating set of cardinal finite objects. The main theorem is then a generalization obtained as an application of this result internally in the localic reflection of an arbitrary topos: a topos is generated by cardinal finite objects if and only if it is separated, locally decidable, and its localic reflection is zero dimensional.
\end{abstract}

\renewcommand{\thefootnote}{\fnsymbol{footnote}} 
\footnotetext{\emph{Keywords.} topos, separated topos, finite object}
\footnotetext{\emph{2010 Mathematics Subject Classification.} 18B25, 03G30.}
\renewcommand{\thefootnote}{\arabic{footnote}} 


\tableofcontents

\section{Introduction}

\blockn{Cardinal finite objects in a topos, i.e. objects that are locally constant and finite, are extremely useful when it comes to do analysis over this topos. By that we mean, for example, studying bundle of Hilbert, Banach or maybe Frechet spaces over it, defined as Hilbert, Banach or Frechet spaces in the internal logic.

For example, in some of our previous works (\cite{henry2014measure}, \cite{henry2015toward}) we have shown that for a boolean topos $\Tcal$ the existence of finite objects in sufficient quantity is equivalent to the fact that $\Tcal$ is separated, and we have used this to prove several interesting results about continuous fields of Hilbert spaces over locally separated boolean toposes.}

\blockn{The goal of this paper is to extend this result of existence of cardinal finite objects to non-boolean toposes, in order to study non-boolean toposes by the same kind of methods that those we have used in the boolean case. The central result is theorem \ref{ThHccase} asserting that, constructively, a hyperconnected separated and locally decidable topos admit a generating family of cardinal finite objects. One can then apply this result internally in the localic reflection of an arbitrary topos and we obtain the main theorem (\ref{MainTh}), which characterizes toposes admitting enough cardinal finite objects as being the toposes which are separated, locally decidable and whose localic reflection is zero dimensional, which generalizes our theorem \ref{ThHccase}.}

\blockn{In a subsequent paper (\cite{henry2015GreenJulg}) we will apply theorem \ref{ThHccase} to obtain a description of the category of Hilbert bundles over a topos which is separated, locally decidable, and whose localic reflection is locally compact, as Hilbert modules over a $C^{*}$-algebra attached to the topos. The idea here is that it is well known how to deal with the localic part for this kind of results and theorem \ref{ThHccase} produces exactly the tools needed to deal with the hyperconnected part of a separated topos. In the same spirit, the two main ingredients that were missing to extend the results of \cite{henry2015toward} to non-boolean topos are the existence of finite objects and a good notion of normal functor between category of Hilbert bundle over a topos. The first is provided by theorem \ref{ThHccase} and the second by the notion of continuous functors between pre-complete $C^{*}$-categories introduced in \cite{henry2015GreenJulg}, hence one can now probably try to prove a reconstruction theorem as in \cite{henry2015toward} for topos which are locally decidable and locally ``separated, with a completely regular localic reflection'' using the exact same method as in the boolean case.}

\blockn{Section \ref{sec_prelim} contains some preliminaries on the framework (in \ref{subsec_framework}) and some general definition of constructive mathematics and topos theory (in \ref{subsec_topos}). Subsection \ref{sec_prelim_coh} recall the definition of coherent locales, their relation to distributive lattices and proves a descent properties for coherent locales.

Section \ref{sec_statement} contains the statement of the main theorem, as well as the proof of the easy direction of the equivalence and the proof of the theorem for localic toposes (which is relatively easy).

Section \ref{sec_hyperco} is in some sense the core of the paper: it contains the proof of the main theorem in the case of a hyperconnected topos (stated as theorem \ref{ThHccase}). It is the most difficult step in the proof of the main theorem. 

Section \ref{sec_endproof} explains how the main theorem for an arbitrary topos can be obtained from the theorem for localic toposes (proposition \ref{localiccase}) and (the internal interpretation of) the theorem for hyperconnected toposes (theorem \ref{ThHccase}) using the hyperconnected-localic fatorization. It concludes the proof of the main theorem \ref{MainTh}.
 }

\section{Preliminaries}
\label{sec_prelim}
\renewcommand{\thesubsubsection}{\arabic{section}.\arabic{subsection}.\arabic{subsubsection}}

\subsection{The logical framework}
\label{subsec_framework}

\blockn{We refer to \cite{maclane1992sheaves} or \cite{borceux3} for an introduction to topos theory and locale theory, and to \cite{sketches} for the more advance topos theoretic notions. We will in particular make an intensive use of the internal logic of toposes (the Kripke-Joyal semantics) as introduced in these references.}

\blockn{The general framework is the internal logic of an elementary topos $\Scal$, called the base topos, endowed with a natural object\footnote{The natural number object is only here for convenience, it is possible to obtain the same results without it.}. In particular, we will use intuitionist logic everywhere. Except for this base topos $\Scal$, the word topos will always mean ``Grothendieck topos'', or more precisely\footnote{The results of \cite[B3]{sketches} show that a topos bounded over $\Scal$ is the same as a Grothendieck topos ``in the internal logic''' of $\Scal$ in the sense that it is defined by an internal site.} a topos bounded over $\Scal$. Objects of $\Scal$ will be simply called sets. One could of course take $\Scal$ to be the category of ordinary sets, but we will need at some point to apply results proved in $\Scal$ to other toposes which is why we restrict ourselves to this framework.}

\blockn{We will also need to use some form of unbounded quantification, i.e. quantification over the class of all sets, which is not part of what is usually called internal logic. To do this we will use M.Shulman ``Stack semantics'' as introduced in \cite[section 7]{shulman2010stack}, which is an extension of ordinary internal logic. We are not making any additional assumption on the base topos, in particular we are not assuming that it is autological in the terminology of M.Shulman \cite{shulman2010stack}. It means that we are only allowed to use the bounded (or restricted, or $\Delta_0$ ) separation axiom, and not the full separation axiom.

In more explicit terms, it means that we can indeed quantify over all sets but this form of quantification cannot be use within the definition of a sub-object. For example, the ``set of $x \in X$ such that there exists an object $S$ satisfying something'' does not necessarily exists, but it make sense to ask whether or not there exists a subset of $X$ which corresponds to this definition. }

\blockn{The external translation of these unbounded quantifications are as follow: 

An existential quantification $\exists X, P(X)$, for some proposition $P$ involving an object $X$, means in a topos $\Tcal$ that there is an inhabited object $S \twoheadrightarrow 1 $ of $\Tcal$ and an object ``$X$'' of $\Tcal_{/S}$ which satisfies the property $P$ in the topos $\Tcal_{/S}$. 
A universal quantification $\forall X, P(X)$ is interpreted as the fact that for all slice $\Tcal_{/V}$ of the topos, for all object $A$ of $\Tcal_{/V}$, the object $A$ satisfies the proposition $P$ in the topos $\Tcal_{/V}$.

Of course in both case $P$ can contain itself some unbounded quantification in which case those definitions are used inductively (together with the rest of the internal logic).

The only uses we will have of this form of quantification will be to say that some family of objects is generating (i.e. that all other objects of the topos can be covered by object in the family), or to say that a topos admit a generating set satisfying some properties. }

\blockn{Strictly speaking, and even using the stack semantics, it is not possible to talk about ``Grothendieck toposes'' in the internal logic, because Grothendieck toposes are large categories and the internal logic can only talk about small objects. This is not a major obstruction, and we will avoid this issue by taking the following conventions:
\begin{itemize}
\item By a ``Grothendieck topos'' $\Tcal$ we actually mean a site.
\item An ``object'' of a topos $\Tcal$ is a sheaf over the corresponding site.
\item A (geometric) morphism $f :\Tcal \rightarrow \Ecal$ of Grothendieck toposes is a flat continuous functor from the site corresponding to $\Ecal$ to the category of sheaves over the site corresponding to $\Tcal$, with the obvious notion of natural transformation between them. It is a classical result (see \cite[VII.10]{maclane1992sheaves}) that those are the same as what is usually called geometric morphism. This defines a (large) $2$-category which is entirely described in terms of ``small'' objects.
\item A site of definition for a topos is a site which defines an isomorphic topos (together with the choice of an isomorphism).
\end{itemize} 
With theses conventions we completely avoid to deal with any sort of large category, and we are still able to talk about Grothendieck topos in the internal logic in a transparent way as if we had universes. Moreover, it is true that a topos $\Tcal$ (bounded\footnote{But morphismes between Grothendieck toposes are always bounded.}) over another topos $\Ecal$ is the same things as a topos in the internal logic of $\Ecal$, in the sense that the corresponding $2$-categories are equivalent. }

\subsection{Topos theoretic preliminaries}
\label{subsec_topos}

\block{A set indexed collection $(X_i)_{i \in I}$ of objects of a topos is said to be generating if it satisfies one of the following equivalent conditions: every object can be covered by a co-product of the $X_i$; the $X_i$ form a family of generators (and/or strong generators, and/or regular generators) in the usual categorical sense; there is a site of definition of the topos such that the $X_i$ are exactly the representable sheaves.}

\block{For a set $X$ one says that:

\begin{itemize}
\item $X$ is (Kuratowski) \emph{finite} if $\exists n, \exists x_1, \dots ,x_n \in X$ such that $\forall x \in X, \exists i, x=x_i$.

\item $X$ is \emph{cardinal finite} if there exists $n$ such that $X \simeq \{1,\dots ,n \}$

\item $X$ is \emph{decidable} if for all $x,y \in X$ one has $x=y$ or $x \neq y$.

\item $X$ is \emph{inhabited} if $\exists x \in X$.

\item A subset $Y \subset X$ is complemented if for all $x \in X$, $x \in Y$ or $x \notin Y$.

\end{itemize}

An object of a topos is said to be finite, cardinal finite, decidable or inhabited if it satisfies those properties internally in the logic of the topos. We remind the reader that those are local properties, for example, for a sheaf over a topological space, being inhabited does not mean to have a global section but only to have local sections around any point of the topological space, and being cardinal finite does not mean to be a constant sheaf, but only locally constant: a cardinal finite sheaf is the same as a finite etale cover, and if the space is connected and locally connected they corresponds to finite sets endowed with an action of the $\pi_1$.
}

\block{We recall a few simple properties of these notions in intuitionist mathematics:

\begin{itemize}
\item A subset of a (cardinal) finite set does not have to be finite.
\item A quotient of a finite set is finite.
\item A set is cardinal finite if and only if it is finite and decidable.
\item A finite subset of a decidable set is complemented.
\item A complemented subset of a finite set is finite.
\item A finite set is either empty or inhabited (depending on if $``n"$ is zero or not).
\item A set is not inhabited if and only if it is empty but a non empty set has no reasons to be inhabited.
\item For a cardinal finite set $X$ the integer $n$ such that $X$ is isomorphic to $\{1,\dots, n \}$ is unique and is called the cardinal of $X$, we denote it by $|X|$.
\end{itemize}

More on the theory of (Kuratowski) finite objects can be found in \cite[D5.4]{sketches}, this reference also explains how to define finiteness without a natural number object. If one want to work in this framework (i.e. without a natural number object) one has to define cardinal finite as ``decidable and Kuratowski finite", remove all mention of cardinality of a set and replace the induction over cardinal that appears in the paper by induction over the Kuratowski finite sets themselves. We will not do that in the present paper for the sake of simplicity.
}

\block{If $\Tcal$ is a topos and $X$ is an object of $\Tcal$, $\Tcal_{/X}$ denotes the slice topos whose objects are objects of $\Tcal$ endowed with a morphism to $X$ and morphisms are commutative triangle. In the case of the topos of sheaves over a topological space it is the topos of sheaves over the etale space of the sheaf $X$.}

\block{\label{generating_collection}A topos is said to be \emph{locally decidable} if it has a generating set of decidable objects and it is said to be \emph{generated by cardinal finite objects} if it has a generating set of cardinal finite objects. It appears, that, in our framework, this is equivalent to a weaker definition:

\Prop{A (Grothendieck) topos $\Tcal$ is locally decidable (resp. generated by cardinal finite objects) if and only for every object $X$ of $\Tcal$, there is a set $I$, a collection $(D_i)_{i \in I}$ of decidable (resp. cardinal finite) objects of $\Tcal$ and an epimorphism:
\[ \coprod_{i \in I} D_i \twoheadrightarrow X \]
 }
 
\Dem{If there is a generating set of decidable (resp. cardinal finite) objects then one can take $I$ to be the set of all maps from an element of the generating set to $X$ and one obviously obtains that the condition given in the proposition holds.

The other direction relies on the use of collection axiom, and more precisely on the axiom of `` collection of context'' stated in proposition $5.6$ of \cite{shulman2010stack}. This axiom is valid internally in the stack semantics of any elementary topos (in fact any extensive Heyting category) by lemma $7.13$ and proposition $5.6$ of \cite{shulman2010stack}. 

Let $\Tcal$ be a Grothendieck topos $\Tcal$ satisfying the condition stated in the proposition, and let $(A_j)_{j\in J}$ be a set of generators for $\Tcal$. For every $j \in J$ there exists a covering of $A_j$ by a co-product of decidable (resp. cardinal finite) objects, hence by the axiom of collection of contexts mentioned above, there exists a set $S$, an epimorphism $u:S \twoheadrightarrow J$, a family of set $(I_s)_{s \in S}$, a family of decidable (resp. cardinal finite) objects $(D_{i,s})_{i \in I_s}$ and a family of epimorphisms:

\[ w_s : \coprod_{i \in I_s} D_{i,s} \twoheadrightarrow A_{u(s)} \]

Hence the $D_{i,s}$ for $s \in S$ ad $i \in I_s$ form a set of decidable (resp. cardinal finite) generators for $\Tcal$.

}

Our definition of locally decidable is slightly weaker than usual notion which would be ``any object is a quotient of a decidable object'', because a co-product of decidable objects can fail to be decidable if the set of indices of the co-product is not decidable. 

To clarify the discussion of the relation between the two notions, we will say that a topos is quasi-decidable if every object is a quotient of a decidable object, but this notion is not going to be used in the present paper. A topos is always internally (in itself) locally decidable, as internally the terminal object is a decidable generator. But one can see that a topos is internally (in itself) quasi-decidable if and only if it is locally quasi-decidable. The two notions are equivalent if and only if the base topos is internally quasi-decidable i.e. locally quasi-decidable. Also note that despite the name ``locally decidable'', neither local decidability nor quasi-decidability are local properties: Take $\Tcal$ to be a non-boolean locally decidable topos (for example $\Tcal =\sh([0,1])$) and $G$ be a non-decidable group object in $\Tcal$ (for example $G$ be the quotient of $\{+1,-1\}$ whose two global sections are identified on $(0,1]$) then the topos $\Tcal_{G}$ of $G$-object of $\Tcal$ is not locally decidable, but it has a (covering) slice isomorphic to $\Tcal$ hence is ``locally locally decidable''. In our example, $\Tcal_{G}$ is the topos of sheaf over $[0,1]$ endowed with an involution which is identity on all the stalks except possibly the stalk at $0$, and the functor that forget the action of $G$ is the ``$f^*$ part'' of an etale covering of $\Tcal_{G}$ by $\Tcal$.

}

\block{A geometric morphism $f :\Ecal \rightarrow \Tcal$ is said to be localic if internally in $\Tcal$, $\Ecal$ is a topos of sheaves over some locale. We will not distinguishes between a topos localic over $\Scal$ and the corresponding locale.

As explained in \cite[A4.6]{sketches}, localic morphisms are part of a unique factorization system ``hyperconnected/localic". Hyperconnected morphisms will be of great importance in this paper. The hyperconnected/localic factorization of a morphism $f:\Tcal \rightarrow \Scal$ produces a locale called the localic reflection of $\Tcal$ whose open subspaces are the subterminal objects of $\Tcal$, and $f:\Tcal \rightarrow \Scal$ is hyperconnected if and only its localic reflection is trivial, i.e. if and only if $f_*(\Omega_{\Tcal}) \simeq \Omega_{\Ecal}$ where $\Omega$ denotes the sub-object classifier and the isomorphism is induced by the natural comparison map.}

\blockn{When $X$ is a locale (or a localic topos) $\Ocal(X)$ denotes the corresponding frame (the frame of open subspaces or of sub-terminal objects).}

\block{A geometric morphism if said to be \emph{proper} if its localic part is, internaly in the target, a compact locale (i.e. satisfies the usual finite sub-cover property). A geometric morphism $f: \Ecal \rightarrow \Tcal$ is said to be \emph{separated} if the diagonal map $\Ecal \rightarrow \Ecal \times_{\Tcal} \Ecal$ is proper. A topos is said to be \emph{separated} if its canonical geometric morphism to $\Scal$ is separated. These notions have been introduced by I.Moerdijk and C.J.Vermeulen in \cite{moerdijk2000proper}, and a lots of their properties have been proved there. These notions also play an essential role in the present paper.}

\blockn{We will also need the notion of regular locale as defined for example in \cite{sketches} between C1.2.16 and C1.2.17 and the fact that regular locales are separated (see \cite{sketches} C1.2.17). A locale is said to be regular if any open subspace $U \subset X$ can be covered by open subspaces $V \subset U$ such that the closure of $V$ is included in $U$.}

\block{If $f :\Ecal \rightarrow \Fcal$ is any geometric morphism and $\Tcal$ is a topos internally in $\Fcal$, then $f^{\sharp}(\Tcal)$ denotes the topos in $\Ecal$ that is obtained by pulling back a site of definition of $\Tcal$. In terms of topos over $\Fcal$ and $\Ecal$ it corresponds to the pullback $\Ecal \times_{\Scal} \Tcal$ seen as a topos over $\Ecal$. When $\Tcal$ is a locale, $f^{\sharp}(\Tcal)$ is a locale internally in $\Ecal$ (i.e. the map to $\Ecal$ is localic) and has (internally in $\Ecal$) a basis of its topology given by $f^{*}(\Ocal(\Tcal))$.}

\subsection{Spectrums of distributive lattices and descent}
\label{sec_prelim_coh}

\blockn{In this section we briefly review the theory of coherent and locally coherent locales in order to prove the descent property of \ref{coherenceDescend}.}

\block{By a \emph{distributive lattice} we mean a partially ordered set $E$ in which every finite subset has a supremum (denoted by the symbol $\vee$ ) and every pair of elements has an infimum (denoted by $\wedge$) and which satisfies the distributivity law:

\[(a \vee b) \wedge c = (a \wedge c) \vee (b \wedge c). \]

We \emph{do not assume} that our distributive lattices have a top element. A distributive lattice is said to be \emph{unital} if it has a top element. Distributive lattices always have a bottom element given by the empty supremum.
}

\block{Let $E$ be a distributive lattice, the spectrum of $E$ is the locale $\spec E$ defined by:

\[ \Ocal(\spec E) = \{ \text{Ideals of } E \}\]

Where an ideal of $E$ is a subset $I$ of $E$ stable under finite supremum, containing the bottom elements and such that if $a \in I$ and $b \leqslant a$ then $b\in I$.

$\spec E$ classifies the theory of prime filters of $E$, i.e. subsets $J \subset E$ which are inhabited, stable under finite intersection, if $a \in J$ and $a \leqslant b$ then $b \in J$, the bottom element of $E$ is not in $J$, and if $a \vee b \in J$ then $a \in J$ or $b \in J$.

Assuming the axiom of choice, $\spec E$ is indeed the ordinary topological space of prime ideal of $E$, but in constructive mathematics such prime ideals do not need to exists and $\spec E$ can have no points.
}

\block{\Def{A locale is said to be locally coherent if it has a basis of compact open subspaces which is stable under intersection. It is said to be coherent if it is locally coherent and compact.
A morphism $f$ between two locally coherent locales is said to be coherent if $f^{*}$ send compact open subspaces to compact open subspaces.
}}

\block{\Prop{The spectrum of a distributive lattice is a locally coherent locale, whose compact open subspaces correspond to the elements of $E$ (they are the principal ideals).

Conversely, any locally coherent locale is the spectrum of its distributive lattice of compact open subspace.

This induces an anti-equivalence of categories between the category of distributive lattices and locally coherent locales (with coherent morphisms between them). Under this correspondence unital distributive latices correspond to coherent locales.}

\Dem{See \cite[II.3]{johnstone1986stone} for the unital case, whose proof is essentially constructive. The proofs for the non-unital case are exactly the same.}

}

\blockn{In this ``non-unital'' framework, a boolean algebra is\footnote{Using the usual condition that every element has a complement would imply the existence of a top element.} a distributive lattice in which for every pair of elements $a \leqslant b$ there is an element $c=b-a$ such that $a \wedge c $ is the bottom element and $a \vee c =b$ (such an element is unique). One can still check (but we will not need it) that boolean algebras are the same as (non-unital) commutative ring in which every element is a projection. }

\block{\Prop{A locally coherent locale is separated if and only if it is regular, and if and only if it is the spectrum of a boolean algebra.}

\Dem{If $\spec E$ is separated, then for any compact subspace $K$ of $\spec E$ the map from $K$ to $\spec E$ is proper (by \cite[II.2.1.(iv)]{moerdijk2000proper} ) and hence $K$ is closed in $\spec E$. In particular, for any pair of elements $a \leqslant b$ in $E$, the open subspace corresponding to $a$ is also closed and hence it has an open and closed complement. The intersection of this complement with $b$ is a compact open subset $c$ (hence an element of $E$) such that $c \vee a = b$ and $c \wedge a = \emptyset$.

Conversely, the spectrum of a boolean algebra has a basis of closed open subspaces and hence is regular and hence separated.
}
}

\block{\label{specPullbackstable}\Prop{Let $f:\Ecal \rightarrow \Scal$ be any geometric morphism, let $E$ be a distributive lattice in $\Scal$ then there is a canonical isomorphism $f^{\sharp} \spec E  \simeq \spec f^{*} E$.}

\Dem{This follows directly from the description of $\spec E$ in term of classifying space: $f^{\sharp} \spec E$ and $\spec f^{*} E$ classify the same theory: the prime filters of $f^{*}E$.}

}

\blockn{The last proposition of this section relies heavily on descent theory for locales and objects (it is stated for opens and proper surjections, but it actually holds for any descent morphism for locales and we will need it only for hyperconnected geometric morphism). Section C5.1 of \cite{sketches} contains everything that we need about descent theory, we will follow the terminology of this reference (in particular, ``descent" means ``effective descent").}

\block{\label{coherenceDescend}\Prop{Let $f : \Ecal \rightarrow \Scal$ be a geometric morphism which is either an open surjection or a proper surjection. Let $\Lcal$ be a locale in $\Scal$ such that, internally in $\Ecal$, $f^{\sharp}\Lcal$ is locally coherent then $\Lcal$ is also locally coherent in $\Scal$.}

The same result also holds if one replaces ``locally coherent'', by coherent or (locally) coherent and separated either by using exactly the same proof as below (replacing distributive lattices by unital distrbutive lattices or boolean algebras) or by using \cite[C.5.1.7]{sketches}, which among other thing proves similar results for compactness and separation.

\Dem{Any isomorphism between two locally coherent locales is coherent, hence induces an isomorphism between the corresponding distributive lattices. Using this simple observation together with \ref{specPullbackstable}, one see that the canonical descent data on the locally coherent locale $f^{\sharp}\Lcal$ induces a descent data on the corresponding distributive lattice $B$. As open and proper surjections are descent morphisms for objects, this proves that the distributive lattice $B$ (together with its canonical descent data) is of the form $f^{*}B'$ for $B'$ a distributive lattice in $\Scal$, by definition $f^{\sharp}\Lcal \simeq \spec B \simeq \spec f^{*} B' \simeq f^{\sharp}\spec B'$ and all these objects carries a descent data and the isomorphisms are compatible to the descent data, hence (as proper and open surjections are descent morphisms for locales) this descend into an isomorphism between $\Lcal$ and $\spec B'$ which concludes the proof.

}

}

\section{Statement of the main result}
\label{sec_statement}
\renewcommand{\thesubsubsection}{\arabic{section}.\arabic{subsubsection}}

\blockn{The main result of this paper is the following theorem, which holds in constructive mathematics:}

\block{\label{MainTh}\Th{A topos $\Tcal$ is generated by cardinal finite objects if and only if it satisfies all of the following conditions:

\begin{itemize}

\item $\Tcal$ is locally decidable.

\item $\Tcal$ is separated.

\item Every sub-terminal object of $\Tcal$ can be written as a union of complemented sub-terminal objects.

\end{itemize}

}

The necessity of the first condition is clear: if there is a generating set of cardinal finite objects then there is in particular a generating set of decidable objects. The necessity of the last condition is also relatively easy: if $X$ is a cardinal finite object of a topos, the support of $X$ (i.e. the sub-terminal object corresponding to the proposition ``$X$ is inhabited'') is complemented\footnote{Because a a Kuratowski finite set is either empty or inhabited}, but as there is a generating set of cardinal finite objects, any subterminal object is a union of supports of such objects, and hence of complemented sub-objects. Corollary \ref{sepLEM2} below proves the necessity of the second condition. The proof of the theorem will be finished at the very end of the present paper.
}

\block{Let us also observe that a topos satisfying the conclusion of the theorem is a coherent topos if and only it is compact, i.e. if and only if its localic reflection is compact, and that in this case the coherent objects are exactly the cardinal finite objects: this follows from proposition IV.3.4 of \cite{moerdijk2000proper}, and the fact that compact and separated implies strongly compact ($III.2.8$ of \cite{moerdijk2000proper}).}

\block{\Lem{\label{lemmafiniteimpsep}Let $\topos$ be a topos which admit a generating set of cardinal finite objects. Let $p$ and $q$ be two points of $\topos$. Then the locale $I$ of isomorphisms from $p$ to $q$, defined as the pullback:

\[
\begin{tikzcd}[ampersand replacement=\&]
I \arrow{r} \arrow{d} \& * \arrow{d}{p} \\
* \arrow{r}{q} \& \topos \\
\end{tikzcd}
\]

is compact.
}

\Dem{This locale $I$ classifies the theory of isomorphisms from $p$ to $q$. By hypothesis, one can construct a site of definition $\Ccal$ for $\topos$ whose representable objects correspond to cardinal finite objects. As points, $p$ and $q$ can then be seen as flat continuous functors from $\Ccal$ to sets (see \cite[VII.10]{maclane1992sheaves}), and a morphism from $p$ to $q$ can then be described as a collection of maps from $p(c)$ to $q(c)$ for $c \in \Ccal$ satisfying the natural transformation condition. Because of the assumption, all the $p(c)$ and $q(c)$ are cardinal finite sets, hence the theory of collection of maps from $p(c)$ to $q(c)$ is compact by the localic Tychonof theorem\footnote{See for example \cite{coquand2003compact} for a proof of a statement of this theorem general enough for our purpose.} and the condition of being a natural transformation can be written as a series of equality which, because equality is decidable in cardinal finite objects, is going to be an intersection of closed sublocales hence is again closed and is hence compact.

An isomorphism is given by a couple of morphisms which are inverses of each other hence the classifying locale for isomorphisms is a closed sublocale of the locale hom$(p,q) \times$hom$(q,p)$ and hence is also compact, which concludes the proof.
}

}

\block{\label{sepLEM2}We can now prove that separation is a necessary condition in theorem \ref{MainTh}.

\Cor{Let $\topos$ be a topos admitting a generating set of cardinal finite objects then $\topos$ is separated.}

\Dem{The key observation, is that the hypothesis of being generated by a set of cardinal finite objects is pullback stable, in the sense that if $\topos$ satisfies it and $\Ecal$ is another topos then $\topos \times \Ecal$ satisfies it internally in $\Ecal$ because a pullback of a cardinal finite object is cardinal finite and a site for $\topos \times \Ecal$ internally in $\Ecal$ is given by pulling back to $\Ecal$ a site for $\topos$ in sets. Hence as the lemma \ref{lemmafiniteimpsep} has been proven constructively its conclusion will holds internally in $\Ecal$ for all the the pullback of $\topos$. The trick is then as usual to apply this to the universal case: that is internally in $\Ecal = \topos \times \topos$ in which (the pullback of) $\topos$ has two canonical points, namely $\pi_1$ and $\pi_2$ and the locale of isomorphisms between them is exactly the diagonal embeddings of $\topos$ in $\topos \times \topos$, hence one gets that this map is proper and that $\topos$ is separated. 
}
}

\blockn{We can also immediately prove the localic version of theorem \ref{MainTh}. In this case, only the last hypothesis of the theorem is relevant: localic toposes are always locally decidable, and the separation condition is implied by the third condition (because it implies regularity). }

\block{\label{localiccase}\Prop{Let $\Lcal$ be a localic topos, then $\Lcal$ is generated by cardinal finite objects if and only if it has a basis of complemented open subspaces.}

\Dem{We already observed the necessity of this condition right after the statement of theorem \ref{MainTh}. Conversely, if $\Lcal$ has a basis of complemented open subspaces, then those complemented open subspaces form a generating set of cardinal finite objects, which concludes the proof. }
}

\section{Proof in the hyperconnected case}
\label{sec_hyperco}

\blockn{In \cite{moerdijk2000proper}, I.Moerdijk and C.J.Vermulen proved that every hyperconnected separated topos with a point is equivalent to the topos of sets endowed with a continuous action of the (compact) localic group $G$ of automorphisms of the point. It might be tempting to think that this representation theorem immediately implies that those toposes are generated be cardinal finite objects: indeed, as $G$ is compact the orbit of any point should be a finite sets. This is not true. Indeed, in the absence of the law of excluded middle, the locale $G$ does not have to be open\footnote{Or overt, or locally positive, this mean that the canonical map from $G$ to the point is an open map.} and this introduces all sort of complications: for example, if $X$ is a set endowed with an action of $G$ and $x$ is a point of $X$ then the orbit of $x$ under the action of $G$ does not have to be a subset of $X$ (it is a sublocale of the discrete locale $X$), and there is no\footnote{In fact, the existence of orbits would implies that the topos is atomic and this would implies in turn that $G$ is open because of corollary C3.5.14 of \cite{sketches} and under this assumption the simpler proof used in the boolean case can be applied.} ``minimal'' $G$-stable subset containing $x$, and hence no natural candidate for a cardinal finite sub-object included in $X$ and containing $x$. The first half of this section would be about constructing such cardinal finite objects in the topos of $G$-sets.}

\blockn{Note that if we are only interested in theorem \ref{ThHccase} in the framework of classical mathematics one can avoid all these additional difficulties and find proofs\footnote{This is done for example in section 3 of \cite{henry2014measure}.} that are considerably shorter and simpler than the one given here, but this will not be enough for proving theorem \ref{MainTh} or for the others applications we have in mind.}

\block{\label{lemcompactor}\Lem{Let $K$ be a compact locale and $p$ its canonical map to the point, and let $U$ be a proposition (i.e. a sub-set of the singleton). Assume that $p^{*}(U) = K$ then one has $U$ or $``K=\emptyset"$.}

Also note that without any compactness assumption on $K$, one has that the proposition $``K=\emptyset"$ is equal to $p_*(\emptyset)$. Indeed, for any proposition $V$, $V \leqslant p_*(\emptyset)$ if and only if $p^{*}(V)=\emptyset$ if and only if $V \Rightarrow ( K=\emptyset )$.

\Dem{As $p^{*}$ commute with arbitrary supremums, one has:
\[ p^{*}(U) = p^{*}\left( \bigvee_{* \in U} \{*\} \right) = \bigvee_{* \in U} K \]

hence if $p^{*}(U)=K$ one can deduce by compactness of $K$ that there is a finite subset $V$ of $U$ such that the union indexed by $V$ is already equal to $K$. A finite set is either empty or inhabited, in the first case this proves that $K = \emptyset$ in the second case that $U$ is inhabited, hence that the proposition $U$ holds.
}

}

\block{\label{emptycompactproduct}\Lem{Let $K_1$ and $K_2$ be two compact locales such that $K_1 \times K_2 = \emptyset$ then $K_1 = \emptyset$ or $K_2 = \emptyset$.}

\Dem{Consider the following pullback square:
\[
\begin{tikzcd}[ampersand replacement=\&,scale=2]
K_1 \times K_2 \arrow{r}{\pi_1} \arrow{d}{\pi_2} \&  K_1 \arrow{d}{p_1} \\
K_2 \arrow{r}{p_2} \&  \{*\}\\
\end{tikzcd}
\]
By proposition C.3.2.6 of \cite{sketches} one has the following commutative Beck-Chevalley square:

\[
\begin{tikzcd}[ampersand replacement=\&,scale=2]
\Ocal(K_2) \arrow{r}{(p_2)_*} \arrow{d}{\pi_2^{*}} \&  \Ocal(\{*\})  \arrow{d}{p_1^{*}}\\
\Ocal(K_1 \times K_2) \arrow{r}{(\pi_1)_*}  \&  \Ocal(K_1).\\
\end{tikzcd}
\]

Because $K_1 \times K_2 = \emptyset$, $(\pi_1)_*$ is constant equal to the top element of $\Ocal(K_1)$ hence $(p_1)^{*}(p_2)_* \emptyset = K_1$ hence by lemma \ref{lemcompactor} one has $K_1 = \emptyset$ or $(p_2)_*(\emptyset)$ but $(p_2)_*(\emptyset)$ is the proposition $K_2=\emptyset$ hence this concludes the proof.
}
}

\block{\label{finitnessForCompactGrp}\Prop{Let $G$ be a compact localic group acting on a decidable set $X$. Let $x \in X$ then there exists a finite subset $F \subset X$ stable under the action of $G$ such that $x \in F$. }

\Dem{Let $G$ and $X$ be as in the proposition.

For any $x \in X$ let $\mu_x$ be the map of ``action on $x$'' $G \rightarrow X$ defined on generalised elements by $g \mapsto g.x$.

By compactness of $G$ for any $x \in X$ there is a (cardinal) finite subset $x \in S \subset X$ which contains the image of $\mu_x$.

Fix an $x \in X$ let $x \in S \subset X$ be such a finite set and let $S \subset S' \subset X$ be a finite set which contain the image of $\mu_s$ for all $s \in S$.

We will now prove by induction on $|S'|+|S|$ that given a pairs $(S,S')$ of finite subset of $X$ such that $Gx \subset S$ and $GS \subset S'$ there exists a finite subset $F$ such that $x \in F \subset S'$ and $F$ is stable by the action of $G$.

If $|S'| = |S|$ it means that $S = S'$ and hence $S$ is already stable under the action of $G$ (this cover in particular the ``initialisation case'')

Assume that $|S'| > |S|$, and let $v \in S'$ which is not in $S$ (let us recall that as a finite subset of a decidable set $S$ is complemented).

For any pair $a,b \in X$ let $G_{a,b}=(\mu_a)^{-1}(\{ b \})$ i.e., in term of generalized elements $G_{a,b} = \{g \in G | g.a = b \}$. As $X$ is decidable $\{ b \}$ is both open and closed as a sublocale of the discret locale $X$ and hence $G_{a,b}$ is open and closed in $G$, in particular it is compact. The multiplication map of $G$ induces a morphism $G_{b,c} \times G_{a,b} \rightarrow G_{a,c}$ hence by lemma \ref{emptycompactproduct}, if $G_{a,c} = \emptyset$ then for any $b$ either $G_{a,b} = \emptyset$ or $G_{b,c} = \emptyset$.

Coming back to our $v \in S'$, because $v \notin S$, $G_{x,v}=\emptyset$ hence for each $s \in S$ either $G_{x,s} =\emptyset$ or $G_{s,v} = \emptyset$, by finiteness of $S$ there is two possible case: 

\begin{itemize}
\item For all $s \in S$ one has $G_{s,v}=\emptyset$:

In this case we will prove that the pair $(S,S'-\{v\})$ still satisfies our induction hypothesis, this will implies (by induction) that there exists a finite subset $F \subset S'-\{v\} \subset S'$ containing $x$ and stable under $G$ and concludes the proof.
 
For any $s\in S$, $\mu_s^{-1}(\{v\}) =G_{s,v} = \emptyset$ hence $\mu_{s}^{-1}(S'-\{v\})=\mu_{s}^{-1}(S')=G$ hence $Gs \subset S'-\{v\}$, which proves that $(S,S'-\{v\})$ indeed satisfies our hypothesis.

\item There is an $s \in S$ such that $G_{x,s}=\emptyset$.

In this case we will prove that the pair $(S-\{s\},S')$ satisfies the induction hypothesis, this will implies (by induction) that there exists a finite subset $F \subset S'$ containing $x$ and stable under $G$ and concludes the proof.

As $ 1 \in G_{x,x}$, $x$ is different from $s$ hence $x \in S-\{s\}$ moreover $\mu_x^{-1}(\{s\})= \emptyset$ hence $Gx \subset S-\{s\}$. Finally $G(S-\{s\}) \subset GS \subset S'$ so this concludes the proof.

\end{itemize}

}

}

\block{\label{Isspectrum1}\Cor{Let $\Tcal$ be a separated hyperconnected topos with a point, and $X \in \Tcal$ be decidable object, then the localic reflection of $\Tcal_{/X}$ is a locally coherent separated locale.}

\Dem{Let $p$ be a point of $\Tcal$. By theorem II.3.1 of \cite{moerdijk2000proper}, $\Tcal$ is equivalent to the topos of sets endowed with an action of $G$, the compact localic group of automorphisms of $p$, and any object $X$ of $\Tcal$ is associated to the set $p^{*}(X)$ endowed with its canonical action of $G$. In particular, if $X$ is a decidable object of $\Tcal$, $p^{*}(X)$ is a decidable set endowed with an action of a compact localic group, hence one can apply proposition \ref{finitnessForCompactGrp} and $p^{*}X$ can be covered by finite $G$-stable subsets. This means that $X$ can be covered by sub-object which are internally finite. If $U \subset X$ is such a sub-object, then $U$ is complemented in $X$ (because it is finite in a decidable set), $\Tcal_{/U} \rightarrow \Tcal$ is compact because $U$ is internally finite (see \cite{moerdijk2000proper}[I.1.4]), $\Tcal \rightarrow *$ is proper because $\Tcal$ is hyperconnected hence $\Tcal_{/U}$ is compact (a composite of proper map is proper, \cite{moerdijk2000proper}[I.2.1]) i.e. the frame of sub-objects of $U$ is compact. So $U$ corresponds to a compact complemented open subspace of the localic reflection of $\Tcal_{/X}$.

Finite subobjects in a decidable object are stable under binary intersections (because they are complemented and a complemented sub-object of a finite object is finite) hence the localic reflexion of $\Tcal_{/X}$ has a basis of complemented compact open subspaces stable under binary intersections which concludes the proof.
}
}

\block{\label{unpointedcoherence}\Prop{Let $\Tcal$ be a hyperconnected separated topos and $X$ a decidable object of $\Tcal$. Then the localic reflection of $\Tcal_{/X}$ is locally coherent and separated.}

\Dem{Let $f: \Ecal \rightarrow \Scal$ be an open or proper surjection such that internally in $\Ecal$, $f^{\sharp}(\Tcal)$ has a point. For example, one can take $\Ecal = \Tcal$ which is both an open and a proper surjection because $\Tcal$ is hyperconnected and $\Tcal$ has a point in its internal logic given by the diagonal map $\Ecal = \Tcal \rightarrow  \Tcal \times \Tcal = f^{\sharp} \Tcal$.

The pullback by $f$ preserves all the hypothesis: $f^{\sharp} \Tcal$ is separated by \cite[II.2.2.(i)]{moerdijk2000proper} and hyperconnected by \cite[C2.4.11(ii)]{sketches}, $f^{\sharp}(\Tcal_{/X}) \rightarrow f^{\sharp}\Tcal$ is etale and separated (hence corresponds internally to the slice by a decidable object of $f^{\sharp}\Tcal$ by \cite[II.1.3.(i)]{moerdijk2000proper}) and by \cite[C2.4.12]{sketches} the hyperconnected/localic factorization $\Tcal_{/X} \rightarrow \Lcal \rightarrow \Scal$ becomes a hyperconnected/localic factorization:  $f^{\sharp}\Tcal_{/X} \rightarrow f^{\sharp} \Lcal \rightarrow \Ecal$.
Applying corollary \ref{Isspectrum1} internally in $\Ecal$ gives that $f^{\sharp} \Lcal$ is (internally) locally coherent and separated and hence proposition \ref{coherenceDescend} shows that $\Lcal$ is itself locally coherent and separated, which concludes the proof.

}
}

\blockn{We can finally prove the main result of this section:}

\block{\label{ThHccase}\Th{Let $\Tcal$ be a hyperconnected separated locally decidable topos, then $\Tcal$ admit a generating set of cardinal finite objects.}

In particular, by \cite[Prop. IV.3.4]{moerdijk2000proper} such a topos is coherent and its coherent objects are exactly the internally cardinal finite objects.

\Dem{$\Tcal$ has a generating set of decidable objects. We will prove that any decidable object of $\Tcal$ is covered by its finite sub-objects and this will conclude the proof.

Let $X$ be a decidable sub-object of $\Tcal$. By proposition \ref{unpointedcoherence}, the localic reflection of $\Tcal_{/X}$ is locally coherent, hence $X$ is covered by subobjects $U \subset X$ such that $\Tcal_{/U}$ is compact. But as $\Tcal$ is separated, this proves (by \cite[II.2.1.(iv)]{moerdijk2000proper}) that the morphism $\Tcal_{/U} \rightarrow \Tcal$ is proper, hence $U$ is internally Kuratowski finite, but as a sub-object of $X$ it is decidable, hence it is cardinal finite, which concludes the proof.
}

}

\section{Proof in the general case}
\label{sec_endproof}

\blockn{We can now finish the proof of theorem \ref{MainTh}. We just need to explain a little more what it means externally to be internally locally decidable or generated by cardinal finite objects.}

\block{\label{Prop_characLocdec}\Prop{Let $f:\Tcal \rightarrow \Ecal$ be a geometric morphisms. Then the following conditions are equivalent:

\begin{enumerate}

\item Internally in $\Ecal$, $\Tcal$ is locally decidable (resp. generated by cardinal finite objects).

\item For any object $X$ of $\Tcal$, there exists an object $I \in \Ecal$, a decidable (resp. cardinal finite) object $D$ in $\Tcal_{/f^*I}$ and an epimorphism $D \twoheadrightarrow X$.

\end{enumerate}
}

\Dem{ We first assume $(1)$. Let $X$ be any object of $\Tcal$, then internally in $\Ecal$, $X$ admit a covering by the co-product of a family $(D_i)_{i \in I}$ of decidable (resp. cardinal finite) objects. Externally, it means that there is an inhabited object $S$ of $\Ecal$, an object $I$ of $\Ecal_{/S}$, a decidable (resp. cardinal finite) object $D$ of $\Tcal_{f^* I}$, and an epimorphism $D \rightarrow X \times f^*S$ over $f^*S$. As $S$ is inhabited, $f^*S$ is also inhabited hence the projection map $X \times f^* S \rightarrow X$ is an epimorphism and hence one has a covering of $X$ by an object $D$ which is decidable (resp. cardinal finite) over an object $f^*I$ which concludes the proof.

\bigskip

We now assume $(2.)$. Let $(B_j)_{j \in J}$ be a set of generators for $\Tcal$ (over the base, not over $\Ecal$).

For each $j \in J$ there exists an object $T$ of $\Ecal$, a decidable (resp. cardinal finite) object $D$ of $\Tcal_{/f^* T}$ and an epimorphism $D \rightarrow B_j$. Hence, by the collection axiom already discussed in the proof of \ref{generating_collection} there exists a set $J'$ with a surjection $s:J' \rightarrow J$, a collection $(T_{j})_{j \in J'}$ of objects of $\Ecal$, a collection $(D_j)_{j \in J'}$ of decidable (resp. cardinal finite) objects of $\Tcal_{/T_{j}}$ and a $J'$-indexed collection of epimorphisms $D_j \rightarrow B_{s(j)}$.

In particular, the $D_j$ for $j \in J'$ form a generating family of $\Tcal$. Let:
 \[ \Dgo:= \coprod_{j \in J' } D_j, \]
 \[ T:= \coprod_{j \in J'} T_j. \]
 
We claim that $\Dgo$ is, internally in $\Ecal$, a $T$-indexed generating family of decidable (resp. cardinal finite) objects, this claim immediately implies that $(1)$ holds hence concludes the proof.

The map from $\Dgo$ to $f^* T$ given by the co-product of the maps from $D_j$ to $f^* T_j$ indeed makes $\Dgo$ into a $T$-indexed collection of decidable (resp. cardinal finite) objects of $\Tcal$. We have already mentioned that the family $(D_j)_{j \in J'}$ is externally generating so it is also internally generating (seen as a $p^*J'$ indexed family of objects) but internally, in $\Ecal$, for each $j \in J'$, one has:

\[D_j \simeq \coprod_{x \in T_j} \Dgo_{\iota_j(x)} \]

where $\iota_j$ denotes the natural injection from $T_j$ to $T$, hence $\Dgo$ is also internally a generating family and this concludes the proof.

}
}

\block{\label{HCfactislocallydec}\Cor{Let $\Tcal$ be a locally decidable topos, and let $f : \Tcal \rightarrow \Ecal$ be any geometric morphism then, internally in $\Ecal$,  $\Tcal$ is again locally decidable.}

\Dem{This follows immediately from the characterization $(2)$ in proposition \ref{Prop_characLocdec}: If $\Tcal$ is itself locally decidable then any object $X$ of $\Tcal$ admit a covering by a family of decidable objects, hence a covering by an object $Y$ which is decidable over a constant object $p^* I $, where $p$ is the morphism $\Tcal \rightarrow \Scal$, but $p^* I$ is in particular of the form $f^* S$ for $S$ a (constant) object of $\Ecal$ so this concludes the proof. }
}

\block{\label{compoOffinitelyGen}\Cor{Let $f:\Tcal \rightarrow \Ecal$ be any geometric morphism, assume that $\Ecal$ is generated by cardinal finite objects and that, internally in $\Ecal$, $\Tcal$ is generated by cardinal finite objects. Then $\Tcal$ is generated by cardinal finite objects.}

\Dem{It also follows rather directly from characterization $(2)$ of \ref{Prop_characLocdec}. Let $p$ and $p'$ be the morphism $p:\Ecal \rightarrow \Scal$ and $p':\Tcal \rightarrow \Scal$. Let $X$ be an object of $\Tcal$, as $\Tcal$ is, internally in $\Ecal$, generated by cardinal finite objects, there exists a covering of $X$ by an object $Y$ in $\Tcal$ which is cardinal finite over an object $f^*V$ for $V$ an object of $\Ecal$. As $\Ecal$ is also generated by cardinal finite object, $V$ admit a covering by an object $W$ which is cardinal finite over a constant object $p^* J$.

Let $Y' := f^* W \times_{f^*V} Y$. As $W \rightarrow V$ is a covering, $Y' \rightarrow Y$ is a covering and hence $Y'$ is a covering of $X$. As $Y$ is cardinal finite in $\Tcal_{/f^* V}$ , $Y'$ is also cardinal finite in $\Tcal_{/f^*W}$, and as $f^* W$ is cardinal finite in $\Tcal_{/p'^* J}$, one has that $Y'$ is cardinal finite in $\Tcal_{/p'^* J}$ and this concludes the proof.  }

}

\block{At this point, theorem \ref{MainTh} is proved: Let $\Tcal$ be a topos satisfying the three conditions of theorem \ref{MainTh}, and let $\Tcal \rightarrow \Lcal$ be its hyperconnected localic factorization. $\Lcal$ satisfies the condition of proposition \ref{localiccase} hence is generated by cardinal finite objects, internally in $\Lcal$, $\Tcal$ is separated because of proposition \cite[II.2.3]{moerdijk2000proper} and locally decidable because of lemma \ref{HCfactislocallydec}, hence it is generated by cardinal finite objects by theorem \ref{ThHccase}, and hence proposition \ref{compoOffinitelyGen} allows to conclude that $\Tcal$ is generated by cardinal finite objects.}

\bibliography{Biblio}{}

\begin{thebibliography}{10}

\bibitem{borceux3}
F.~Borceux.
\newblock {\em Handbook of {C}ategorical {A}lgebra 3: {S}heaf {T}heory},
  volume~3.
\newblock {Cambridge University Press}, 1994.

\bibitem{coquand2003compact}
Thierry Coquand.
\newblock Compact spaces and distributive lattices.
\newblock {\em Journal of Pure and Applied Algebra}, 184(1):1--6, 2003.

\bibitem{henry2014measure}
Simon Henry.
\newblock Measure theory over boolean toposes.
\newblock {\em {arXiv preprint arXiv:1411.1605}}, 2014.

\bibitem{henry2015GreenJulg}
Simon Henry.
\newblock Complete ${C}^{*}$-categories and a topos theoretic {G}reen-{J}ulg
  theorem.
\newblock {\em arXiv preprint arXiv:1512.03290}, 2015.

\bibitem{henry2015toward}
Simon Henry.
\newblock Toward a non-commutative {G}elfand duality: {B}oolean locally
  separated toposes and monoidal monotone complete ${C}^*$-categories.
\newblock {\em {arXiv preprint arXiv:1501.07045}}, 2015.

\bibitem{johnstone1986stone}
Peter~T. Johnstone.
\newblock {\em Stone spaces}.
\newblock Cambridge University Press, 1986.

\bibitem{sketches}
P.T. Johnstone.
\newblock {\em Sketches of an elephant: a topos theory compendium}.
\newblock Clarendon Press, 2002.

\bibitem{maclane1992sheaves}
S.~MacLane and I.~Moerdijk.
\newblock {\em Sheaves in {G}eometry and {L}ogic: {A} {F}irst {I}ntroduction to
  {T}opos {T}heory}.
\newblock Springer, 1992.

\bibitem{moerdijk2000proper}
Ieke Moerdijk and Jacob Johan~Caspar Vermeulen.
\newblock {\em Proper maps of toposes}, volume 705.
\newblock AMS Bookstore, 2000.

\bibitem{shulman2010stack}
Michael~A. Shulman.
\newblock Stack semantics and the comparison of material and structural set
  theories.
\newblock {\em arXiv preprint arXiv:1004.3802}, 2010.

\end{thebibliography}
\bibliographystyle{plain}

\end{document}